\documentstyle[amscd]{amsart}
\input{bull-art}
\newtheorem{thm}{Theorem}
\newtheorem{lem}[thm]{Lemma}
\theoremstyle{definition}

\newcommand{\thmref}[1]{Theorem~\ref{#1}}
\newcommand{\lemref}[1]{Lemma~\ref{#1}}

\newcommand{\gp}{Gelfand pair}
\newcommand{\C}{\hbox{${\Bbb C}$}}
\newcommand{\R}{\hbox{${\Bbb R}$}}
\newcommand{\T}{\hbox{${\Bbb T}$}}

\newcommand{\g}{\hbox{${\frak g}$}}
\newcommand{\k}{\hbox {${\frak k}$}}
\newcommand{\co}{\hbox {${\cal O}$}}
\newcommand{\Vr}{\hbox{$V_{\Bbb R}$}}
\newcommand{\obc}{\hbox{{(OC)}}}
\newcommand{\fobc}{\hbox{{(FOC)}}}

\begin{document}
\def\currentvolume{31}
\def\currentissue{2}
\def\currentyear{1994}
\def\currentmonth{October}
\def\copyrightyear{1994}
\def\currentpages{185-190}

\title
{The Moment Map for a Multiplicity Free Action}

\author[C. Benson]{Chal Benson}
\address[C. Benson and G. Ratcliff]{Department of 
Mathematics and Computer
Science\\
University of Missouri--St. Louis\\
St. Louis, Missouri 63121}
\thanks{All authors were supported in part by the National 
Science
Foundation}
\email[C. Benson]{benson@@arch.umsl.edu\newline\indent
{\defaultfont{\it E-mail address}, G. Ratcliff: 
}ratcliff@@arch.umsl.edu}
\author[J. Jenkins]{Joe Jenkins}
\address[J. Jenkins]{Department of Mathematics and 
Statistics\\
State University of New York at Albany, Albany, New York 
12222}
\email{jwj71@@math.albany.edu}

\author[R. Lipsman]{Ronald L. Lipsman}
\address[R. Lipsman]{Department of Mathematics\\
University of Maryland\\
College Park, Maryland 20742}
\email{rll@@clio.umd.edu}

\author[G. Ratcliff]{Gail Ratcliff}

\date{January 13, 1994}



\maketitle

\begin{abstract}
Let $K$ be a compact connected Lie group acting unitarily 
on a 
finite-dimensional complex vector space $V$. One calls 
this a
{\em multiplicity-free}
action whenever the $K$-isotypic components of $\C[V]$ are 
$K$-irreducible.
We have shown that this is the case if and only if the 
moment map
$\tau:V\rightarrow\k^*$ for the action is finite-to-one on 
$K$-orbits.
This is equivalent to a result concerning \gp s associated 
with Heisenberg
groups that is motivated by the Orbit Method. Further 
details of this work
will be published elsewhere.
\end{abstract}


\section{Introduction}

\subsection{Gelfand pairs}
Suppose that $G$ is a Lie group and that $K$ is a compact 
subgroup of $G$.
One says that $K\subset G$ is a {\em Gelfand pair} when 
the algebra $L^1(G//K)$
of integrable $K$-bi-invariant functions on $G$ is 
commutative under convolution.
This note concerns a special class of \gp s that are 
associated with Heisenberg groups. 
Let $V$ be a finite-dimensional
complex vector space with Hermitian structure 
$\langle\cdot,\cdot\rangle := \left(\cdot,\cdot\right) + 
i\omega(\cdot,\cdot)$,
and form the associated Heisenberg group $H_V := 
V\times\R$ with product
$(z,t)(z',t') := (z+z',t+t'-\frac 12\omega(z,z'))$. Let 
$K$ be a compact connected Lie
group acting on $V$ via some unitary representation 
$K\times V\rightarrow V$,
$(k,z)\mapsto kz$. We obtain an action of $K$ on $H_V$ by 
automorphisms
$k\cdot(z,t) := (kz,t)$ and form the semidirect product 
$G:=K\ltimes H_V$.
In this setting, $L^1(G//K)$ is naturally identified with 
the algebra $L^1_K(H_V)$
of integrable $K$-invariant functions on the Heisenberg 
group. We will say that
$(K,H_V)$ is a \gp\ when this algebra is commutative.

A fundamental result due to I. M. Gelfand asserts that 
$K\subset G$ is a \gp\
if and only if each irreducible unitary representation 
$\pi$ of $G$
has at most a one-dimensional space of $K$-fixed vectors 
\cite{Gelf}.
The {\em $K$-spherical} representations are those that 
have nonzero
$K$-fixed vectors.
By using the Mackey machine to describe the unitary dual
of $K\ltimes H_V$, one obtains an interesting 
specialization of this criterion
to pairs of the form $K\subset K\ltimes H_V$. Namely, 
$(K,H_V)$ is a \gp\ if and
only if the action of $K$ on $V$ is {\em 
multiplicity-free} \cite{BJR1,Car}.
The latter condition means that each irreducible 
representation of $K$ occurs at
most once in the representation $k\cdot p(z) := 
p(k^{-1}z)$ of $K$ on $\C[V]$.
Multiplicity-free actions are of considerable interest in 
their own right,
but are usually studied in the setting of reductive
complex algebraic groups. One can pass to this setting by 
complexifying $K$.
The {\em irreducible} multiplicity-free actions of 
connected reductive
complex algebraic groups were classified by V. Kac 
\cite{Ka}. This yields a
complete classification of the \gp s $(K,H_V)$ where $K$ 
acts 
irreducibly on $V$ \cite{BJR1}.

The literature contains various results relating 
multiplicities of Lie
group representations to geometric properties of coadjoint 
orbits.
See for example \cite{CG3,GS2,GS3,Heck,Lip2}. 
These results motivate a conjecture that provides a 
geometric
formulation of Gelfand's criterion for \gp s $(K,H_V)$.
 Let \k, \g\ denote the Lie algebras of $K$, $G=K\ltimes 
H_V$,
and let $\k^*$, $\g^*$ be their duals. The {\em Orbit 
Method} (also called 
``Geometric Quantization'') yields  one-to-one 
correspondences between the
unitary duals $\widehat{K}$, $\widehat{G}$ of $K$, $G$ and 
the integral
orbits in $\k^*$, $\g^*$ for the coadjoint actions of $K$, 
$G$.
(See \cite{Heck} and \cite{Lip1} respectively.) Let 
$\pi:\g^*\rightarrow\k^*$
be the restriction map and $\k^\perp:=\pi^{-1}(\{0\})$. If 
$\xi\in\g^*$ and
$\co^G_\xi := \text{Ad}^*(G)\xi$ is an integral orbit 
corresponding to 
$\rho_\xi\in\widehat{G}$, then the multiplicity of the 
trivial representation
of $K$ in $\rho_\xi|_K$ (i.e. the dimension of the space 
of $K$-fixed vectors)
should be related to the number of $K$-orbits in the 
intersection 
$\k^\perp\cap\co_\xi^G$.

\medskip
\noindent {\bf Orbit Conjecture}. {\em Let $K$ be a 
compact connected Lie group 
acting unitarily on $V$. 
Then $(K,H_V)$ is a \gp\ if and only if}
\begin{equation}
\text{for every}\ \xi\in\k^\perp,\ \co_{\xi}^G\cap\k^\perp
\ \text{is\ a\ single}\ K\text{-orbit.}\tag{OC}
\end{equation}

Although this conjecture seems rather bold, it does hold 
in many
interesting cases and we have found no counterexamples to 
date.
Moreover, we have recently proved a result that is only 
slightly weaker
than the Orbit Conjecture.

\begin{thm}\label{thm1}
Let $K$ be a compact connected Lie group 
acting unitarily on $V$.  
Then $(K,H_V)$ is a \gp\ if and only if
\begin{equation}
\text{for\ every}\ \xi\in\k^\perp,\ \co_{\xi}^G\cap\k^\perp
\ \text{is\ a\ finite\ union\ of}\ K\text{-orbits.}
\tag{FOC}\end{equation}
\end{thm}

A closely related result, due to V. Guillemin and S. 
Sternberg, can be found in
\cite{GS3}. They prove that if $G$ is {\em compact}, then 
$K\subset G$ is a \gp\ if and
only if condition $\fobc$\ holds generically. 
\thmref{thm1} does not, however, follow
from the results in \cite{GS3} since in our case, $G$ is 
not compact. We will
outline our proof of this theorem below. Details of this 
work will appear
elsewhere \cite{BJLR}.

\subsection{The moment map}
One can write points $\xi\in\g^*$ as 
$\xi=(\alpha,z_\circ,\lambda)$ where
$\alpha\in\k^*$, $z_\circ\in V$, $\lambda\in\R$, and 
$\xi(A,z,t) = \alpha(A)+\omega(z_\circ,z)+\lambda t$ for 
$A\in\k$, $z\in V$,
and $t\in\R$. One can verify that condition \obc\ always 
holds for ``nongeneric''
orbits of the form $\co_{(\alpha,z_\circ,0)}^G$. Moreover, 
one computes
that for $\lambda\ne 0$, 
$\co_{(\alpha,0,\lambda)}^G\cap\k^\perp
= \left\{(0,z,\lambda)~|~\tau(z)\in\co_{-2\lambda\alpha}^K%
\right\}$
where $\tau:V\rightarrow\k^*$ is the (unnormalized) {\em 
moment map} defined
by
$$\tau(z)(A) := \omega(z,Az)$$
for $z\in V$, $A\in\k$. Here $(A,z)\mapsto Az$ denotes the 
derived action of
\k\ on $V$.
A key property of $\tau$ is $K$-equivariance: 
$\tau(kz)=\text{Ad}^*(k)\tau(z)$.
In particular, $\tau$ maps $K$-orbits in $V$ to 
$\text{Ad}^*(K)$-orbits in $\k^*$.
The above remarks lead to the following reformulations of 
conditions \obc\
and \fobc.
\begin{itemize}
\item Condition \obc\ holds if and only if $\tau$ is 
one-to-one on $K$-orbits.
\item Condition \fobc\ holds if and only if $\tau$ is 
finite-to-one on $K$-orbits.
\end{itemize}
We can now remove the Heisenberg group from the picture, 
restating
\thmref{thm1} as an equivalent result that provides a 
geometric criterion
for multiplicity-free actions.

\begin{thm}\label{thm2}
The action of $K$ on $V$ is multiplicity-free if and only 
if the moment map
$\tau:V\rightarrow\k^*$ is finite-to-one on $K$-orbits.
\end{thm}

The Orbit Conjecture is equivalent to the assertion that 
$K\times V\rightarrow V$
is a multiplicity-free action if and only if $\tau$ is 
one-to-one on $K$-orbits.
In view of \thmref{thm2}, this conjecture would be proved 
if one could show that
$\tau$ is necessarily one-to-one on $K$-orbits whenever it 
is finite-to-one
on $K$-orbits. Although we have been unable to show this, 
we have proved that
if $\tau$ is finite-to-one on $K$-orbits, then
$\tau$ is necessarily {\em uniformly} finite-to-one on 
$K$-orbits. This means that
for some constant $M$, $\tau^{-1}(\co_\alpha^K)$ contains 
at most $M$ $K$-orbits for
all $\alpha\in\k^*$. Thus, the finiteness result needed to 
complete the proof of the
Orbit Conjecture would provide a geometric counterpart to 
the following lemma,
which has a folklore status in the field.

\begin{lem}\label{lem1}
If the  multiplicities of the representations of $K$ 
occurring in $\C[V]$ are 
finite and bounded, then $K\times V\rightarrow V$ is a 
multiplicity-free
action.
\end{lem}

\subsection{Capelli actions}
Let $I[\k^*]$ denote the $\text{Ad}^*(K)$-invariant 
polynomials on $\k^*$ and
$\C[\Vr]^K$ denote the $K$-invariant polynomials on the 
underlying real
vector space \Vr\ for $V$. Since $\tau$ is 
$K$-equivariant, we obtain
an associated map $\tau^*:I[\k^*]\rightarrow \C[\Vr]^K$, 
$\tau^*(p)(z):=p(\tau(z))$. Condition \obc\ is equivalent 
to the assertion
that $\tau^*(I[\k^*])$ separates $K$-orbits in $V$. Since 
$K$ is compact, one
knows that the full algebra $\C[\Vr]^K$ of invariant 
polynomials always separates
$K$-orbits. Thus we see that condition \obc\ certainly 
holds whenever
$\tau^*(I[\k^*])=\C[\Vr]^K$. We say that $K\times 
V\rightarrow V$ is a
{\em Capelli action} in this case. \thmref{thm2} shows 
that Capelli actions
are always multiplicity free. Moreover, these provide an 
interesting class
of Gelfand pairs that satisfy the Orbit Conjecture.

In \cite{HU} R. Howe and T. Umeda consider the ``abstract 
Capelli problem''
for actions of reductive complex algebraic groups. 
$K\times V\rightarrow V$ 
is a Capelli action if and only if the abstract Capelli 
problem for its
complexification $K_{\text c}\times V\rightarrow V$ has an 
affirmative answer.
This means that the derived action of the universal 
enveloping algebra
${\cal U}(\k_c)$ on $V$ yields a surjective map 
$d\iota:{\scr{Z\/U}}(\k_c)\rightarrow {\scr{PD}}(V)^K$. 
Here ${\scr{Z\/U}}(\k_c)$
denotes the center of ${\cal U}(\k_c)$, and 
${\scr{PD}}(V)^K$ is the algebra
of $K$-invariant polynomial coefficient differential 
operators on $\C[V]$.
This relation between Capelli actions and the abstract 
Capelli 
problem follows immediately from the commutative diagram 
(1) below.

The compact forms
of many of Kac's irreducible multiplicity-free actions 
\cite{Ka} are
Capelli actions. From \cite{HU} one obtains a complete 
list of the
irreducible Capelli actions (up to equivalence):
$$\begin{array}{llllll}
U(n) & \text{on}\ \C^n &
U(n) & \text{on}\ S^2(\C^n)  \text{\ ($n\ge 2$)} \\
U(n) & \text{on}\ \Lambda^2(\C^n)  \text{\ ($n\ge 2$)} & 
\T\times SO(n,\R)  & \text{on}\ \C^n  \text{\ ($n\ge 3$)} \\
U(n)\times U(m) & \text{on}\ \C^n\otimes\C^m &
\T\times \text{Sp}(n) & \text{on}\ \C^{2n} \\
U(2)\times \text{Sp}(n) & \text{on}\ \C^2\otimes\C^{2n} &
U(n)\times \text{Sp}(4) & \text{on}\ \C^n\otimes \C^8  
\text{\ ($n\ge 4$)} \\
\T\times \text{Spin}(7)  & \text{on}\ \C^8 &
\T\times \text{Spin}(10)  & \text{on}\ \C^{16} \\
\T\times G_2 & \text{on}\ \C^7
\end{array}$$

\subsection{Spectrum of $\protect\text{Ind}_K^G(1_K)$}
Suppose that $(K,H_V)$ is a \gp. The $K$-spherical 
representations of $G$
are precisely the {\em spectrum} of irreducible unitary
representations weakly contained in the quasi-regular 
representation
$\text{Ind}_K^G(1_K)$
of $G$ on $L^2(G/K)$. The Orbit Method suggests that these 
representations
should correspond to the integral coadjoint orbits in 
$\g^*$ that
meet $\k^\perp$. Equivalently, one might expect that 
$\sigma\in\widehat{K}$ should occur in $\C[V]$ if and only 
if the 
corresponding coadjoint orbit in $\k^*$ belongs to the 
image of the
moment map $\tau:V\rightarrow\k^*$. It is not hard to 
verify this for
the standard action of the unitary group $U(V)$ (with 
respect to
$\langle\cdot,\cdot\rangle$) on $V$. This fact together 
with functoriality
of the moment map and a theorem of G. Heckman \cite{Heck} 
concerning projections of
coadjoint orbits establishes the following result.

\begin{thm}
If $\sigma\in\widehat{K}$ occurs in $\C[V]$, then the 
corresponding
coadjoint orbit in $\k^*$ is contained in $\tau(V)$.
\end{thm}

Although the converse does  hold for some interesting 
examples,
it is not true in general. $\tau(V)$ can contain integral 
coadjoint orbits
that do not correspond to representations of $K$ occurring 
in $\C[V]$.
This can happen even when $K\times V\rightarrow V$ is a 
Capelli action.
The action of $U(n)$ on the symmetric 2-tensors 
$S^2(\C^n)$ provides an
example of this phenomenon. It is known that the 
representations of $U(n)$ appearing
in the (multiplicity-free) decomposition of 
$\C[S^2(\C^n)]$ are parameterized
by Young's diagrams with all rows of even length. On the 
other hand, 
we have verified that
the integral coadjoint orbits in the image of $\tau$ yield 
all representations
parameterized by Young's diagrams of even size (total 
number of cells).
Thus, the image of $\tau$ produces infinitely many 
representations 
of $U(n)$ that do not appear in $\C[S^2(\C^n)]$.

\section{Outline of proof of \thmref{thm2}}

A key ingredient in our proof is the following commutative 
diagram
relating $\tau^*$ to the action of the center of the 
complexified
enveloping algebra of \k\ on $\C[V]$ by polynomial 
coefficient differential 
operators.
\begin{equation}\label{comdiag}
\begin{CD}
I[\k^*]  @>{\tau^*}>> \C[\Vr]^K \\
@VV{\lambda}V         @VV{\Gamma}V \\
{\cal{Z\/U}}(\k_c) @>{d\iota}>> {\cal{PD}}(V)^{K}
\end{CD}
\end{equation}
Here $\lambda$ is the symmetrization map (up to scalars), 
and $\Gamma$ 
takes elements of $V^*$ to
multiplication operators and elements of $\overline{V}^*$ to
differentiation operators.
 The horizontal maps are algebra maps, and the vertical
maps are vector space isomorphisms.

Suppose that $K\times V\rightarrow V$ is a 
multiplicity-free action.
It follows that ${\cal{PD}}(V)^{K}$ is abelian \cite{HU} 
and hence,
 by a result of F. Knop \cite{Knop}, ${\cal{PD}}(V)^{K}$
is finitely generated as a module over
$d\iota({\cal{Z\/U}}(\k_c))$. Although $\Gamma$ is {\em 
not} an algebra
map, we do know that the terms of highest order in 
$\Gamma(pq)$ and
$\Gamma(p)\Gamma(q)$ agree for any $p,q\in\C[\Vr]^K$. An 
induction argument
on the total degree of operators in ${\cal{PD}}(V)^{K}$ 
shows that $\C[\Vr]^K$
is a finitely generated $\tau^*(I[\k^*])$-module.

Let $\gamma_1, \gamma_2, \dots, \gamma_\ell$ be algebra 
generators for
$\C[\Vr]^K$. Each $\gamma_j$ is algebraic over the 
fraction field
of $\tau^*(I[\k^*])$. Suppose that $\gamma_j$ satisfies a 
polynomial
equation of degree $d_j$ with coefficients 
$a_{i,j}\in\tau^*(I[\k^*])$. If
$z_\circ\in V$ and $z\in\tau^{-1}(\tau(Kz_\circ))$, then 
$\gamma_j(z)$
is a root of the degree $d_j$ polynomial with
coefficients $a_{i,j}(z_\circ)$.
Thus, $\gamma_1\times\gamma_2\times\cdots\times\gamma_%
\ell:V\rightarrow\C^\ell$
assumes at most $M:=d_1d_2\cdots d_\ell$ distinct values on 
$\tau^{-1}(\tau(Kz_\circ))$. As the level sets of 
$\gamma_1,\gamma_2,\dots,\gamma_\ell$ are the $K$-orbits 
in $V$, this
shows that $\tau$ is (uniformly) finite-to-one on 
$K$-orbits.

Next suppose that $\tau$ is finite-to-one on $K$-orbits 
and choose
\R-valued generators $q_1, \dots, q_m$ for 
$\tau^*(I[\k^*])$.
For simplicity, we assume here that $q_1, \dots,
q_m$ are algebraically
independent over \R.
The level set of $q_1\times\cdots\times q_m$ through 
$z_\circ\in V$ is 
$\tau^{-1}(\tau(Kz_\circ))$, which is a finite union of 
$K$-orbits.

Let $z_\circ$ be a (``generic'') point in $V$ where 
$D(q_1\times\cdots\times q_m)(z_\circ)$ has maximal rank. 
In fact, one can argue that this maximal rank is precisely 
$m$.
Some open neighborhood of $z_\circ$ is foliated by 
$K$-orbits
of codimension $m$.
The orthogonal complement of the tangent space to $Kz_\circ$
at $z_\circ$ is a dimension $m$ affine subspace $A$ of 
\Vr\ that meets
all $K$-orbits through points $z$ sufficiently close to 
$z_\circ$.
Let $\tilde{q}_j := q_j|_A$. $\tilde{q}_1, \dots, 
\tilde{q}_m$ are
algebraically independent over \R\ since a polynomial 
relation between
the $\tilde{q}_j$'s would imply (by $K$-invariance) a 
polynomial
relation between the $q_j$'s in a neighborhood of $z_\circ$.
Now let $p$ be any \R-valued $K$-invariant polynomial on 
\Vr\ and
consider $\tilde{p}:=p|_A$. By dimensional considerations, 
$\tilde{p}$
is algebraic over $\R\left(\tilde{q}_1, \dots, 
\tilde{q}_m\right)$,
and it follows by $K$-invariance that
 $p$ is algebraic over the fraction field of 
$\tau^*(I[\k^*])$.
We conclude that $\C[\Vr]^K$ is a finitely generated 
$\tau^*(I[\k^*])$-module and hence that ${\cal{PD}}(V)^{K}$
is finitely generated as a 
$d\iota({\cal{Z\/U}}(\k_c))$-module.

Let $\sigma\in\widehat{K}$ and ${\cal H}_\sigma$ be the 
space of
highest weight vectors for $\sigma$ in $\C[V]$. The proof of
Proposition 7.1 in \cite{HU} shows that ${\cal{PD}}(V)^{K}$
must act irreducibly on ${\cal H}_\sigma$ when 
${\cal{PD}}(V)^{K}$
is a finitely generated 
$d\iota({\cal{Z\/U}}(\k_c))$-module. One argues
that if $\left\{D_1, \dots, D_\ell\right\}$ are generators 
for
${\cal{PD}}(V)^{K}$ as a 
$d\iota({\cal{Z\/U}}(\k_c))$-module and 
$p\in{\cal H}_\sigma$, then 
${\cal H}_\sigma\subset \text{Span}\left\{D_1p, \dots, 
D_\ell p\right\}$.
Thus the multiplicity of $\sigma$ in $\C[V]$ is at most 
$\ell$ for
all $\sigma\in\widehat{K}$. \lemref{lem1} now implies that
$K\times V\rightarrow V$ is a multiplicity-free action.

\section*{Acknowledgment}
 The authors thank Roger Howe for his many helpful
comments concerning this work.

\bibliographystyle{plain}

\begin{thebibliography}{10}

\bibitem{BJLR}
C.~Benson, J.~Jenkins, R.~Lipsman, and G.~Ratcliff.
 {\em A geometric criterion for {Gelfand} pairs associated 
with the
  {Heisenberg} group},
 preprint.

\bibitem{BJR1}
C.~Benson, J.~Jenkins, and G.~Ratcliff,
 {\em On {Gelfand} pairs associated with solvable {Lie} 
groups},
 Trans. Amer. Math. Soc. {\bf 321} (1990), 85--116.

\bibitem{Car}
G.~Carcanno,
 {\em A commutativity condition for algebras of invariant 
functions},
 Boll. Un. Mat. Ital. {\bf 7} (1987), 1091--1105.

\bibitem{CG3}
L.~Corwin and F.~Greenleaf,
 {\em Spectrum and multiplicities for restrictions of 
unitary
  representations in nilpotent {Lie} groups},
 Pacific J. Math. {\bf 135} (1988), 233--267.

\bibitem{Gelf}
I.~M. Gelfand,
 {\em Spherical functions on symmetric spaces},
 Dokl. Akad. Nauk USSR {\bf 70}
 (1950), 5--8; Amer. Math. Soc. Transl. Ser. 2, vol.
37, Amer.
Math. Soc., Providence, RI, 1964, pp. 39--44.

\bibitem{GS2}
V.~Guillemin and S.~Sternberg,
 {\em Geometric quantization and multiplicities of group 
representations},
 Invent. Math. {\bf 67} (1982), 515--538.

\bibitem{GS3}
\bysame,
 {\it Multiplicity free spaces},
 J. Differential Geom. {\bf 19} (1984), 31--56.

\bibitem{Heck}
G.~J. Heckman,
 {\em Projections of orbits and asymptotic behavior of 
multiplicities for
  compact connected {Lie} groups}, Invent. Math. {\bf 67}
(1982), 333--356.

\bibitem{HU}
R.~Howe and T.~Umeda,
 {\em The {Capelli} identity, the double commutant theorem 
and
  multiplicity-free actions},  Math. Ann. {\bf 290}
(1991), 565--619.

\bibitem{Ka}
V.~Kac,
 {\em Some remarks on nilpotent orbits},
 J. Algebra {\bf 64} (1980), 190--213.

\bibitem{Knop}
F.~Knop,
 {\em A {Harish-Chandra} homomorphism for reductive group 
actions},
 preprint.

\bibitem{Lip1}
R.~Lipsman,
 {\em Orbit theory and harmonic analysis on {Lie} groups 
with co-compact
  nilradical},
 J. Math. Pure Appl. {\bf 59} (1980), 337--374.

\bibitem{Lip2}
\bysame,
 {\em Orbital parameters for induced and restricted 
representations},
 Trans. Amer. Math. Soc. {\bf 313} (1989), 433--473.

\end{thebibliography}

\end{document}